\documentclass[12pt]{article}
\usepackage[cp1251]{inputenc}
\usepackage{amssymb,latexsym,amsmath}
\usepackage{amsfonts}
\usepackage[russian,english]{babel}
\usepackage{indentfirst}
\usepackage{graphicx,graphics,verbatim}
\usepackage{cmap}
\newcommand{\lr}[1]{\left(#1\right)}
\newtheorem{theorem}{Theorem}
\newtheorem{definition}{Definition}
\newtheorem{corollary}{Corollary}
\begin{document}

\title{On the correctness of finite--rank approximations by series of shifted Gaussians}

\date{}

\maketitle

\author{
\begin{center}
\textbf{S.~M. Sitnik}\\
Belgorod State University, Belgorod, Russia\\
sitnik@bsu.edu.ru\\
\textbf{A.~S. Timashov}\\
Belgorod State University, Belgorod, Russia\\
loaderrus@gmail.com\\
\textbf{S.~N. Ushakov}\\
Voronezh State University, Voronezh,  Russia\\
ushakowww@yandex.ru
\end{center}
}


\begin{abstract}

In this paper we consider interpolation problem connected with series by integer shifts of Gaussians.
Known approaches for these problems met numerical difficulties. Due to it another
method is considered based on finite--rank  approximations by linear systems.
The main result for this approach is to establish correctness of the finite--rank linear system under consideration. And the main result of the paper is to prove correctness of the finite--rank linear system
approximation. For that an explicit formula for the main determinant of the linear system is derived to demonstrate that it is non--zero.
\end{abstract}

\newpage

\section{Introduction and problem statement}

For a long period the main instrument for approximations was based on expansions by complete orthogonal systems. But nowadays in different fields of mathematics and applications occurred more and more problems solutions of which needs  expansions by incomplete, overdetermined or non-orthogonal systems. Such problems are widely considered  for electric or optic signals, filtration, holography, tomography and medicine.  We have such examples of   overdetermined systems as frames, and non-orthogonal as wavelets, Gabor systems and coherent states, Rvachev's systems and so on.

Consider a problem of approximation for arbitrary function defined by its values at integer points.
And stress once more that functions considered are really arbitrary as we work only with function values at integer points.
Let consider approximations by special type series in integer shifts of Gaussians (quadratic exponents with parameters). This system is incomplete in standard spaces but all the same very effective. For the history of this class of approximations, basic results and multiple applications see \cite{S1}--\cite{S3}.

Exactly we will study the next main

\textit{Interpolation problem:} consider arbitrary function $f(x)$ defined on the line $x \in \mathbb{R}$, and a parameter $\sigma>0$, which in probabilistic applications means the standard deviation. We seek for an interpolating function  $\tilde{f}(x)$ also defined on the line $x \in \mathbb{R}$ in the form of the series by integer shifts of Gaussian
\begin{equation}
\label{eq1}
\tilde{f}\left( x \right) \sim \sum\limits_{k=-\infty }^\infty {f_{k}e^{-\frac{\left( x-k \right)^{2}}{2\sigma^{2}}}}
\end{equation}
which coincides with  $f(x)$  at integer points
\begin{equation}
\label{eq2}
f\left( m \right)=\tilde{f}\left( m \right), \   m \in \mathbb{Z}.
\end{equation}

Some approaches to this problem were studied. In \cite{S1} the solution is derived with using special functions, exactly Jacobi theta--functions. But as it was demonstrated in \cite{S2}--\cite{S3} in spite of theoretical value such approach has not computational perspectives as it involves division on very small denominators which leads to large numerical errors. Another approach based on discreet Fourier transform (DFT) is effective but for very small range of parameters \cite{S3}. To overcome these difficulties in  \cite{S4} another approach was considered based on a reduction to finite--rank linear system solution. A serious  obstacle for this last approach was that it was not grounded as a correctness of this method was not established. In this paper we eliminate this gap by proving that a finite--rank linear system is correctly solved. We prove that in a direct and clear way by calculating explicitly the main system determinant and show that it is non--zero. This result is a necessary theoretical ground for a development of numerical methods and algorithms and also rid of difficulties with special functions and DFT.

\section{Reducing to an infinite system of linear equations}

Let introduce a convenient notation

\begin{equation}
\label{eq3}
e \left( \sigma ,x,k \right)=e^{-\frac{{(x-k)}^{2}}{2\sigma^{2}}}
\end{equation}

Solving of the problem formulated is reduced to determination of a sequence of unknown coefficients
$f_{k}$ from (\ref{eq1}). For this following the standard way for interpolation problems we need to find nod functions for every interpolation integer nods $x=m,~ m \in \mathbb{Z}$.  For our case we have to find a single
\textit{basic node function} for the node  $x=0$ which we will seek in the form

\begin{equation}
\label{eq4}
G\left( \sigma, x \right)=\sum\limits_{k=-\infty }^\infty g_{k}\  e(\sigma ,x,k).
\end{equation}

We shortly consider steps to derive an infinite linear system for finding coefficients of the basic node function \eqref{eq4}, cf. \cite{S2}--\cite{S3}.

It follows from (\ref{eq2}) that this basic node function satisfies for all
 $m \in \mathbb{Z}$:
\begin{equation}
\label{eq5}
G \left( \sigma, m \right)=\sum\limits_{k=-\infty }^\infty g_{k}\  e(\sigma, m, k) =\sigma_{m0},
\end{equation}
where $\sigma_{m0}$ is the Kroneker symbol
$$
\sigma_{m0}=\left\{ {\begin{array}{l}
	1,\, m=0 \\
	0,\, m\ne 0. \\
	\end{array}} \right.
$$
If such a function $G\left( \sigma ,x \right)$ obeying a condition (\ref{eq5}) is found then it is easy to find a formal solution of the problem considered. In fact the function
$$
G_{l} \left( \sigma, x \right)= G\left( \sigma, x-l \right)
$$
is the node function for the node $x=l$ as for all $m$
$$
G_{l}\left( \sigma, m \right)= G\left( \sigma, m-l \right)=\sigma_{ml}.
$$
Then one solution of the interpolating problem considered will be a function
\begin{equation}
\label{eq6}
\tilde{f}\left( x \right)=\sum\limits_{l=-\infty }^\infty f\left( l \right) G_{l}\left( \sigma, x \right),
\end{equation}
because for $x=m$ in the sum (\ref{eq6}) only one term is non--zero:
$$
f\left( m \right)G_{m}\left( \sigma, m \right) = f\left( m \right)\cdot
1=f\left( m \right).
$$

To move from solution representation (\ref{eq6}) to the one we seek in the form (\ref{eq1}) do a necessary substitution. By it we derive using (\ref{eq4}):

$$
\tilde{f}\left( x \right)=\sum\limits_{l=-\infty }^\infty f\left( l \right) G_{l}\left( \sigma, x \right)=
$$
$$
=\sum\limits_{l=-\infty }^\infty f \left( l \right) G \left( \sigma, x-l \right) =\sum\limits_{l=-\infty}^\infty f\left( l \right) \sum\limits_{l=-\infty }^\infty g_{k}\ e\left( \sigma, x-l, k \right)
$$

Let introduce a new summation index $j=l+k$ instead of $l=j-k$ and formally change the order of summation. We derive
\begin{equation}
\label{eq7}
\begin{array}{c}
\tilde{f} (x) = \sum\limits_{j=-\infty }^\infty \sum\limits_{k=-\infty }^\infty f  (j-k) g_{k}\ e (\sigma, x-j ,k) = \\
=\sum\limits_{l=-\infty }^\infty \left\{ \sum\limits_{l=-\infty }^\infty f\left( j-k \right)g_{k} \right\}\ e\left( \sigma, x-j, k \right)=  \sum\limits_{j=-\infty }^\infty f_{j}\ e\left( \sigma, x, j \right),
\end{array}
\end{equation}
where coefficients to find are represented in the form (after index change $j \rightleftarrows k$) to agree with
(\ref{eq1})
\begin{equation}
\label{eq8}
f_{k}=\sum\limits_{j=-\infty }^\infty f\left( k-j \right)g_{j},
\end{equation}
and $f(m)$ are values of the given function at integer points, and $g_{j}$ are coefficients of the basic node function series (\ref{eq4}).

Now transform the system of equations:
$$
\sum\limits_{k=-\infty }^\infty g_{k} e \left( \sigma, m, k \right)=\sigma_{m0},~m\in \mathbb{Z}.
$$
For that let introduce a new variable $q=e^{-\frac{1}{2\sigma^{2}}}$. And after that we derive
\begin{equation}
\label{eq9}
\sum\limits_{k=-\infty }^\infty g_{k} q^{{(m-k)}^{2}}=\sigma_{m0},~	m \in \mathbb{Z}.
\end{equation}

The next important step is to consider and study finite--rank cuts of this infinite system of equations (\ref{eq9}).

\section{Reducing to a finite--rank  system of linear equations}

As it was considered before the key moment in the interpolation approach is to define the node function.
Now we pass to study  finite--rank approximations of the initial problem. These approximations are naturally
finite--rank cuts of  infinite system of equations (\ref{eq9}). This direct approach was studied in \cite{S4}. Being in some way restrictive this approach allows to avoid difficulties from other approaches in  \cite{S1}--\cite{S3}. And the most important is that our method allows to construct effective numerical calculations based on it, cf. \cite{S4}.

So let us seek for node function (\ref{eq4}) $G\left( \sigma, x \right)$ approximations of the form $H\left( \sigma, x \right)$ as finite sums
\begin{equation}
\label{eq10}
H\left( n,  x, \sigma \right)=\sum\limits_{k=-n}^n d_{k} \cdot
q^{(x-k)^{2} }, \  q=\exp \left( -\frac{1}{2\sigma^{2}} \right), \  0<q<1,
\end{equation}
after that an infinite system (\ref{eq9}) is reduced to a finite one, and number of equations may be greater than a number of unknown quantities.

\begin{equation}
\label{eq11}
H(n, m,j, \sigma )=\delta_{0j}, \  j=-m,\dots ,0,\dots m, \  m\ge n.
\end{equation}

The system of linear equations which follows from conditions (\ref{eq10})--(\ref{eq11}) there are $2m+1$ equations and $d_{k}$, $-n\le k\le n$ unknown numbers.

Transform the system following from (\ref{eq10}) -- (\ref{eq11}) in the matrix form
\begin{equation}
\label{eq12}
A\cdot \overline d=\overline y,
\end{equation}
with matrix and vector elements
$$
a_{ij}=q^{(i-j)^{2}},~ y_{j}=\delta_{0j}, ~ i=-n, \dots, 0, \dots, n, ~ j=-m, \dots, 0, \dots, m.
$$
For considered for $m=n$ coefficients $d_{k}$ of approximate node function we denote as $H(n, x, \sigma )$, and for $m>n$ as $H(n, m, x, \sigma )$.

Let introduce some notation. Vandermond determinant $n\times n$ we denote as $W \lr{x_1, \dots, x_n}$, this determinant without line $l$ and column $k$ as $W_{l,\, k} \lr{x_1, \dots, x_n}$. For illustrating further considerations let use $n=5$ as examples of arbitrary dimensions, e.g.

$$
W \lr{x_1, x_2, x_3, x_4,  x_5} =
\left|
\begin{array}{ccccc}
1 & x_1 &  x_1^2 &  x_1^3 &  x_1^4 \\
1 & x_2 &  x_2^2 &  x_2^3 &  x_2^4 \\
1 & x_3 &  x_3^2 &  x_3^3 &  x_3^4 \\
1 & x_4 &  x_4^2 &  x_4^3 &  x_4^4 \\
1 & x_5 &  x_5^2 &  x_5^3 &  x_5^4
\end{array} \right|,
$$
$$
W_{3,\, 2} \lr{x_1, x_2, x_3, x_4,  x_5} =
\left|
\begin{array}{ccccc}
1  &  x_1^2 &  x_1^3 &  x_1^4 \\
1  &  x_2^2 &  x_2^3 &  x_2^4 \\
1  &  x_4^2 &  x_4^3 &  x_4^4 \\
1  &  x_5^2 &  x_5^3 &  x_5^4
\end{array} \right|.
$$

It is known that
\begin{equation}
\label{eq13}
W_{l,\, k} \lr{x_1, \dots, x_n} =  \sum x_{\alpha_1} x_{\alpha_2} \dots x_{\alpha_{n-k}}  \cdot \prod\limits_{n \geq i > j \geq 1,~ i \neq l, j \neq l} \lr{x_i - x_j},
\end{equation}
with summation by all sets of $n-k$ numbers $\alpha_1, \alpha_2, \dots, \alpha_{n-k}$  from $1, 2, \dots, n$.

We have also to mention that considered problems are deeply connected with the third Jacobi Theta--function \cite{S5}
$$
\vartheta_{3} \left( z,q \right)=1+2\sum\limits_{k=1}^\infty q^{k^{2}}\cos {\left( 2kz \right),}
$$
and Jacobi--Poisson transformation formula
\begin{equation}
\label{eq14}
\sum\limits_{k=-\infty }^{+\infty } e^{-a\left( t+\pi k \right)^{2}} \, =\frac{1}{\sqrt {\pi a} }\, \, \sum\limits_{k=-\infty }^{+\infty } e^{-\frac{k^{2}}{a}} \, e^{i2kt}.
\end{equation}
For these connections see \cite{S1}--\cite{S3}.

\section{A main case then number of equations equals to a number of unknown coefficients}

\subsection{A correct solvability and formula for the determinant for the $n\times n$ system}

Now the main object of our study is a finite system of linear equations (\ref{eq12}) with square  matrices. We use notation introduced in (\ref{eq10})--(\ref{eq12}).

Let prove our main result of correctness of system solution for square matrix. From it follows that the system is always uniquely solved.

\begin{theorem}
  \label{t1}
	The matrix $A$ for $m=n$ is non--degenerate and its determinant equals to
		\begin{equation}
		\label{eq15}
		|A|=  q^{\frac{2n (n+1)(2n+1)}{3}} \cdot W \lr{q^{-2 n}, \dots 1, \dots,  q^{2 n}}.
	\end{equation}
\end{theorem}

The proof.

$$
|A| =
\left|
\begin{array}{ccccc}
1 & q &  q^4 &  q^9 &  q^{16} \\
q & 1 & q & q^4 & q^9 \\
q^4 & q & 1 & q & q^4 \\
q^9 & q^4 & q & 1 & q \\
q^{16} & q^9 &   q^4 & q & 1
\end{array} \right|
$$
Elements of this determinant is possible to factorize
$$
a_{i j} = q^{(i-j)^2} = q^{{i^2}} \cdot q^{-2ij} \cdot q^{j^2}.
$$
Consequently it is possible to factor out from $i$--th line a value $ q^{{i^2}}$ and from $j$--th column --- $q^{j^2}$. Repeat this operation for all lines and columns
$$
|A| = q^4 \cdot q \cdot 1 \cdot q \cdot  q^4  \cdot \left|
\begin{array}{ccccc}
q^{-4} & q^{-3} &  1 &  q^5 &  q^{12} \\
1 & q^{-1}& 1 & q^3 & q^8 \\
q^4 & q & 1 & q & q^4 \\
q^8 & q^3 & 1 & q^{-1} & 1 \\
q^{12} & q^5 &   1 & q^{-3} & q^{-4}
\end{array} \right| =
q^{20}  \cdot
\left|
\begin{array}{ccccc}
q^{-8} & q^{-4} &  1 &  q^4 &  q^8 \\
q^{-4} & q^{-2}& 1 & q^2 & q^4 \\
1 & 1 & 1 & 1 & 1 \\
q^4 & q^2 & 1 & q^{-2} &  q^{-4}\\
q^8 & q^4 &   1 & q^{-4} & q^{-8}
\end{array} \right|.
$$
Now elements of intermediate determinant are $q^{-2ij}$, further factor out from $i$--th lines multiplicators $q^{2 n i }$, $i = -2 n, \dots, 2 n$ with a unit product:
$$
q^{20}  \cdot
\left|
\begin{array}{ccccc}
q^{-8} & q^{-4} &  1 &  q^4 &  q^8 \\
q^{-4} & q^{-2}& 1 & q^2 & q^4 \\
1 & 1 & 1 & 1 & 1 \\
q^4 & q^2 & 1 & q^{-2} &  q^{-4}\\
q^8 & q^4 &   1 & q^{-4} & q^{-8}
\end{array} \right| =
q^{20}  \cdot
\left|
\begin{array}{ccccc}
1 & q^{4} &  q^8 &  q^{12} &  q^{16} \\
1 & q^{2} &  q^4 &  q^{6} &  q^{8} \\
1 & 1 & 1 & 1 & 1 \\
1 & q^{-2} & q^{-4} & q^{-6} &  q^{-8}\\
1 & q^{-4} &  q^{-8} &  q^{-12} &  q^{-16}
\end{array} \right|=
$$
$$
=
q^{20}  \cdot
\left|
\begin{array}{ccccc}
1 & q^{4} &  \lr{q^{4}}^2 &  \lr{q^{4}}^3 &  \lr{q^{4}}^4 \\
1 & q^{2} &  \lr{q^{2}}^2 &  \lr{q^{2}}^3 &  \lr{q^{2}}^4 \\
1 & 1 & 1 & 1 & 1 \\
1 & q^{-2} &  \lr{q^{-2}}^2 &  \lr{q^{-2}}^3 &  \lr{q^{-2}}^4\\
1 & q^{-4} &  \lr{q^{-4}}^2 &  \lr{q^{-4}}^3 &  \lr{q^{-4}}^4
\end{array} \right|.
$$

In the general case the result is:
$$
\det A = \lr{\prod\limits_{i=-n}^{n} q^{{i^2}}}^2 \cdot \lr{\prod\limits_{j=-n}^{n} q^{{j^2}}}^2  \cdot \prod\limits_{i=-n}^{n}q^{2ni} \cdot W \lr{q^{-2n}, \dots 1, \dots,  q^{2n}} =
$$

$$
=  q^{\frac{2n (n+1)(2n+1)}{3}} \cdot W \lr{q^{-2 n}, \dots 1, \dots,  q^{2 n}}
= q^{\frac{2n (n+1)(2n+1)}{3}} \cdot \prod\limits_{i,j=-n,
	i \neq j}^{n} \lr{q^{-2i}-q^{-2j}}.
$$
Now it is obvious that the main system determinant is non--zero due to inequalities under considered restrictions
$0<q<1$, because of the obvious fact all multipliers in the last formula for the determinant are non--zeros.

Theorem 1 is completely proved.

\textbf{Remark.} The problem to prove correctness of the finite linear system (\ref{eq12}) and so to ground strictly the method of interpolation by shifted Gaussian was first formulated by L.A. Minin and S.M. Sitnik. After rather a time on the basis of numerical computations the formula for the determinant (\ref{eq15}) was found by A.S. Timashov. The strict proof was found by S.N. Ushakov. After that the correctness of the finite linear system (\ref{eq12}) was established and it is a basis for numerical computer methods and algorithms.

\subsection{Further  properties in case of $n\times n$ system: palindromes, symmetry and  dimension reduction}

To formulate further results we define the notion of vector palindrome , this notion is connected with vector components symmetry with respect to its "middle" components.

\begin{definition}
The vector  $\overline x$ in $n$ dimensions is called \textit{a palindrome} if the next relation for its components is fulfilled
$$
\overline{x_i}=\overline{x_{n+1-i}}, \   i= 1,\dots,n.
$$
\end{definition}

\begin{theorem}\label{t2}
Let consider the next linear system of equations
\begin{equation}
	\label{eq16}
	A \cdot \overline x = \overline b,
\end{equation}
where  $A$ is a non--degenerate matrix $n \times n$ and for its elements the next relation is valid
$$
a_{i, j}=a_{n+1-i, n+1-j}\  \forall i,j = 1,\dots,n,
$$
and vector $\overline b$ is a palindrome.  Then the system solution  $\overline x$ is also a palindrome.
\end{theorem}

In fact the condition from theorem 2 means that every matrix $A$ line and column are also palindromes.

The proof.

As the matrix $A$ is non--degenerate then the unique solution $\overline x$ exists. Every $i$--th system line may be written as
$$
\sum\limits_{j=1}^{n} a_{i, j} \cdot x_j = b_i.
$$
Let us prove that a vector $\overline y=\lr{x_{n}, x_{n-1} \dots, x_1 }$ is also a solution to (\ref{eq16}), and from it due to uniqueness property the statement of the theorem 2 follows.

Really for the $i$--th line
$$
\sum\limits_{j=1}^{n} a_{i, j} \cdot y_j = \sum\limits_{j=1}^{n} a_{{n+1-i, n+1-j}} \cdot y_j= \sum\limits_{j=1}^{n} a_{{n+1-i, j}} \cdot x_j = b_{n+1-i}=b_i.
$$

And from it the theorem 2 follows.

Now let us demonstrate the importance of the palindrome property. Using it and symmetry  of the system it is possible to decrease a number of calculations essentially, approximately twice. It leads to more effective and robust numerical algorithms.

\begin{corollary}
For the system (\ref{eq12}) the solution is symmetrical, namely $d_{k}=d_{-k}$.
\end{corollary}

Really let us note that for the system (\ref{eq12}) conditions of the theorem 2 is valid:
$$
a_{ij}=q^{(i-j)^{2}}=q^{(n+1-i-(n+1-j))^{2}}=a_{n+1-i,n+1-j}.
$$
It follows that $d_{k}=d_{-k}$.

This corollary leads to a very important conclusion. Due to it we can reduce a number of equations to solve numerically in fact twice, and it also reduce gaps in values of coefficients. Because of that the numerical difficulties of the problem and computational time are reduced essentially. It permits to solve twice larger systems for the same time as systems without using palindrome symmetry.

Another important fact for the system analysis is that we can find coefficients for approximations of the node function $d_{k}$ explicitly, they are defined by (\ref{eq10}) and are solutions to the system (\ref{eq11}) -- (\ref{eq12}).

\begin{theorem}\label{t3}
	For coefficients   $d_k$ the next formula is valid:
	\begin{equation}
	\label{eq17}
	d_k = (-1)^{k} q^{-k^2} \frac{   W_{k, \, n+1} \lr{q^{-2 n},  \dots, q^0, \dots, q^{2 n}  } }{W \lr{q^{-2 n}, \dots, q^0, \dots, q^{2 n}  }}.
	\end{equation}
\end{theorem}

The proof.

Really by Cramer's rule
$$
d_k = \frac{\Delta_k} {|A|}.
$$
Manipulating with $\Delta_k$ in the same way as in the proof of the theorem 1 we derive
$$
\Delta_k = (-1)^{n+1+k+1+n} q^{-k^2} q^{\frac{2n (n+1)(2n+1)}{3}} W_{k} \lr{q^{- 2 n}, \dots, q^0, \dots, q^{2 n}  },
$$
So
$$
d_k=(-1)^{k} q^{-k^2} \frac{   W_{k} \lr{q^{-2n}, \dots, q^0, \dots, q^{2n}  } }{W \lr{q^{-2n}, \dots, q^0, \dots, q^{2n}  }}.
$$

The theorem 3 is proved.

For example let illustrate the determinant formula \eqref{eq17}  by determinants $5\times 5$ as we do in this paper for clarity and reader's convenience.
$$
\Delta_1 =
\left|
\begin{array}{ccccc}
1 & q &  q^4 &  0&  q^{16} \\
q & 1 & q & 0 & q^9        \\
q^4 & q & 1 & 1 & q^4      \\
q^9 & q^4 & q & 0 & q      \\
q^{16} & q^9 &   q^4 & 0 & 1
\end{array} \right| =
q^{20}  \cdot
\left|
\begin{array}{ccccc}
q^{-8} & q^{-4} &  1 & 0 &  q^8 \\
q^{-4} & q^{-2}& 1 & 0 & q^4 \\
1 & 1 & 1 & q^{-1} & 1 \\
q^4 & q^2 & 1 & 0 &  q^{-4}\\
q^8 & q^4 &   1 & 0 & q^{-8}
\end{array} \right|=
$$
$$
=q^{20} \cdot q^{-1}  \cdot
\left|
\begin{array}{ccccc}
q^{-8} & q^{-4} &  1 &  0&  q^8 \\
q^{-4} & q^{-2}& 1 & 0 & q^4 \\
1 & 1 & 1 & 1 & 1 \\
q^4 & q^2 & 1 & 0 &  q^{-4}\\
q^8 & q^4 &   1 & 0 & q^{-8}
\end{array} \right| =
q^{20}  \cdot q^{-1} \cdot
\left|
\begin{array}{ccccc}
q^{-8} & q^{-4} &  1 &  0&  q^8 \\
q^{-4} & q^{-2}& 1 & 0 & q^4 \\
1 & 1 & 1 & 1 & 1 \\
q^4 & q^2 & 1 & 0 &  q^{-4}\\
q^8 & q^4 &   1 & 0 & q^{-8}
\end{array} \right| =
$$
$$
=
- q^{20} \cdot q^{-1}  \cdot
\left|
\begin{array}{ccccc}
1 & q^{4} &  \lr{q^{4}}^2 &    \lr{q^{4}}^4 \\
1 & q^{2} &  \lr{q^{2}}^2 &   \lr{q^{2}}^4 \\
1 & q^{-2} &  \lr{q^{-2}}^2 &   \lr{q^{-2}}^4\\
1 & q^{-4} &  \lr{q^{-4}}^2 &   \lr{q^{-4}}^4
\end{array} \right|.
$$

And finally let us mention that to  problems considered in this paper the method of transmutations may be applied \cite{S6,S7}. An idea of applying transmutations is to represent a set of Gauss functions as transmuted ones from a set of more standard functions, e.g. linear exponentials. It means that for some transmutation operator $T$ it is valid that
\begin{equation}
\exp(-(x-k)^2) = T (\exp(-(x-k))), k\in \mathbb{Z}.
\end{equation}
So the series in shifted Gaussians may be transformed to series in shifted linear exponents. Such transmutation operators  $T$ may be found explicitly but are rather complicated.

\section{Final conclusions and short paper resume}
In this paper we consider interpolation problem connected with series by integer shifts of Gaussians.
Known approaches for these problems met numerical difficulties. Due to it the direct
method is considered based on finite--rank  approximations by linear systems.
The main result of the paper (theorem 1) is to prove correctness of the finite--rank linear system
approximation. For that an explicit formula for the main determinant of the linear system is derived to demonstrate that it is non--zero. More we consider symmetry properties of the system based on the palindrome notion. It leads to the reduction of the considered system essentially, almost twice in dimension (theorem 2). And at last we derive an explicit formula for the system solution as a ratio of two Vandermond determinants (theorem 3).

\end{document}